\documentclass[11pt]{amsart}
\usepackage{amsmath, amscd, amssymb}
\usepackage{latexsym, amssymb, amsmath, amscd, amsfonts}
\usepackage[mathscr]{eucal}
\usepackage{mathrsfs}
\usepackage{color}
\usepackage{concrete,euler,textcomp}
\usepackage{diagrams, tabmac, url}
\usepackage[T1]{fontenc}

\newtheorem{theorem}{Theorem}[section]

\newtheorem{corollary}[theorem]{Corollary}
\newtheorem{definition}[theorem]{Definition}
\newtheorem{example}[theorem]{Example}

\newtheorem{proposition}[theorem]{Proposition}

\begin{document}
\date{September 8, 2010}
\title{Standard Monomial Theory of RR varieties}
\author{Philip Foth}
\address{CEGEP Champlain - St. Lawrence, Qu\'ebec, G1V 4K2 \ \ Canada  \newline
\ \ and Department of Mathematics\\
The University of Arizona\\
Tucson, AZ 85721 USA}
\email{phfoth@gmail.com}
\author{Sangjib Kim}
\address{
School of Mathematics\\ Tata Institute of Fundamental Research\\ 
Homi Bhabha Road, Mumbai 400005, India}
\email{skim@math.tifr.res.in}

\begin{abstract}
We construct the RR varieties as the fiber products of Bott-Samelson
varieties over Richardson varieties. We study their homogeneous
coordinate rings and standard monomial theory.
\end{abstract}

\maketitle

\section{Introduction}

The main object of our investigation is a fiber product $Z_{w}\times
_{X_{w}^{v}}Z^{v}$, which we call an \textit{RR variety}, over a Richardson
variety $X_{w}^{v}$, of two Bott-Samelson varieties $Z_{w}$ and $Z^{v}$
associated with a Schubert variety $X_{w}$ and an opposite Schubert variety $%
X^{v}$ in the flag variety.%
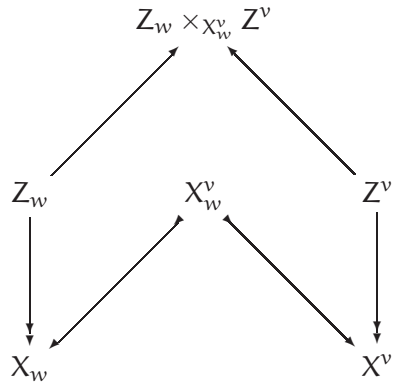
\begin{diagram}
       &        & Z_{w} \times_{X^{v}_{w}} Z^{v} &        &        \\
       &  \ruTo &                                &  \luTo &        \\
Z_{w}  &        &        X^{v}_{w}               &        & Z^{v}  \\
\dOnto &  \ldEmbed &                                &  \rdEmbed & \dOnto \\
X_{w}  &        &                                &        & X^{v}
\end{diagram}We study its standard monomial theory based on an explicit
description of its homogeneous coordinate ring.

\medskip

\section{Schubert and Opposite Schubert varieties}

In this section, we recall definitions and properties of Schubert and
Richardson varieties, and fix our notation.

\subsection{Bruhat-Chevalley order}

Let $G=GL_{n}(\mathbb{C})$ be the general linear group over the complex
numbers $\mathbb{C}$, and $B$ be its Borel subgroup consisting of upper
triangular matrices. We write $T$ for the maximal torus of $G$ consisting of
diagonal matrices. Note that the symmetric group $\mathfrak{S}_{n}$ is the
Weyl group $W=N(T)/T$ of $G$, where $N(T)$ is the normalizer of $T$ in $G$,
and there are finitely many $T$-fixed points $e_{w}$ in $G/B$ labeled by
elements $w$ of $\mathfrak{S}_{n}$. The $B$-orbits $C_{w}=B\cdot e_{w}$ in $%
G/B$ are called the \textit{Schubert cells}. The Zariski closure of $C_{w}$
is called the \textit{Schubert variety} associated with $w$ and denoted by $%
X_{w}$. There is a partial order, called the \textit{Bruhat-Chevalley order}%
, on the elements of $\mathfrak{S}_{n}$:%
\begin{equation*}
w_{1}\geq w_{2}\hbox{ if and only if }X_{w_{1}}\supseteq X_{w_{2}}.
\end{equation*}

The Grassmannian $Gr(d,n)$ of $d$ dimensional subspaces in $\mathbb{C}^{n}$
can be realized as the quotient of the space of $n\times d$ matrices of rank 
$d$ by $GL_{d}(\mathbb{C})$, and then the Schubert varieties in $Gr(d,n)$
can be described explicitly. For $w\in \mathfrak{S}_{n}$, let $i_{k}$\ be
the $k$-th smallest element in $\{w(1),\cdots ,w(d)\}$, then define%
\begin{equation*}
(w(1),\cdots ,w(d))\uparrow =(i_{1},\cdots ,i_{d}),
\end{equation*}%
i.e., the rearrangement of $(w(1),\cdots ,w(d))$ in increasing order. For
the elementary basis $\{e_{i}\}$ of $\mathbb{C}^{n}$, if we let $E_{w}$ be
an $n\times d$ matrix whose $k$-th column represents $e_{i_{k}}$ for $1\leq
k\leq d$, then the Schubert variety $X_{w}$ in $Gr(d,n)$ is the Zariski
closure of its $B$-orbit, $\overline{B\cdot E_{w}}$.

\medskip

The flag variety $G/B$ can be embedded in the product of Grassmannians:%
\begin{equation}
G/B\hookrightarrow Gr(1,n)\times Gr(2,n)\times \cdots \times Gr(n-1,n),
\label{flag-embedding}
\end{equation}%
and the Grassmannian $Gr(d,n)$ can be identified with $G/\widehat{P}_{d}$
for a maximal parabolic subgroup $\widehat{P}_{d}$ of $G$ containing $B$.
Then with respect to the projections $\pi _{d}:G/B\rightarrow G/\widehat{P}%
_{d}$, the Bruhat-Chevalley order can be realized as follows: $w_{1}\geq
w_{2}$ if and only if%
\begin{equation*}
\pi _{d}(X_{w_{1}})\supseteq \pi _{d}(X_{w_{2}})\text{ for }1\leq d\leq n-1
\end{equation*}%
or more explicitly, $w_{1}\geq w_{2}$ if and only if, for $1\leq d\leq n-1$,%
\begin{equation*}
(w_{1}(1),\cdots ,w_{1}(d))\uparrow \succeq (w_{2}(1),\cdots
,w_{2}(d))\uparrow .
\end{equation*}%
where for $a_{1}<\cdots <a_{d}$ and $b_{1}<\cdots <b_{d}$, 
\begin{equation}
(a_{1},\cdots ,a_{d})\succeq (b_{1},\cdots ,b_{d}),  \label{partial_order}
\end{equation}%
if $a_{i}\geq b_{i}$ for all $i$.

With respect to the diagonal embedding (\ref{flag-embedding}), this
condition is compatible with the inclusion order of the Schubert varieties
in each $Gr(d,n)$. We refer to \cite{BL00}\cite{LG01} for further details on
the Schubert varieties and the Bruhat-Chevalley order.

\medskip

\subsection{Flag variety $G/B$}

Recall that for $\mathbf{m}=(m_{1},\cdots ,m_{n-1})\in \mathbb{Z}_{>0}^{n-1}$%
, we have a line bundle $\mathcal{O}_{\mathbf{m}}$ over $G/B$ induced from
the Pl\"{u}cker bundles of the Grassmannians:%
\begin{equation*}
\mathcal{O}_{Gr(1,n)}(m_{1})\otimes \cdots \otimes \mathcal{O}%
_{Gr(n-1,n)}(m_{n-1})
\end{equation*}%
Standard monomial theory lets us describe the section ring of $G/B$:%
\begin{equation*}
\bigoplus_{p\geq 0}H^{0}(G/B,\mathcal{O}_{\mathbf{m}}^{\otimes p})
\end{equation*}%
explicitly in terms of the Pl\"{u}cker coordinates or determinant functions
over the space $M_{n}$ of $n\times n$ complex matrices.

\medskip

To be more precise, let $\mathbb{C}[M_{n}]$ be the coordinate ring of $M_{n}$%
. For $d\leq n$, consider subsets $R=\{r_{1},\cdots ,r_{d}\}$ and $%
C=\{c_{1},\cdots ,c_{d}\}$ of $\{1,\cdots ,n\}$. With $r_{1}<\cdots <r_{d}$
and $c_{1}<\cdots <c_{d}$, we will also write $R$ and $C$ as $(r_{1},\cdots
,r_{d})$ and $(c_{1},\cdots ,c_{d})$.

We let $[R:C]$ or $[r_{1},\cdots ,r_{d}|c_{1},\cdots ,c_{d}]$ denote the map
from $M_{n}$ to $\mathbb{C}$ by assigning to a matrix $X\in M_{n}$ the
determinant of the $d\times d$ minor of $X$ formed by taking rows $R$ and
columns $C$:%
\begin{eqnarray}
\lbrack R:C] &=&[r_{1},\cdots ,r_{d}|c_{1},\cdots ,c_{d}]
\label{RC_determinant} \\
&=&\det \left[ 
\begin{array}{ccc}
x_{r_{1}c_{1}} & \cdots & x_{r_{1}c_{d}} \\ 
\vdots & \ddots & \vdots \\ 
x_{r_{d}c_{1}} & \cdots & x_{r_{d}c_{d}}%
\end{array}%
\right]  \notag
\end{eqnarray}%
We shall identify this determinant with a tableaux obtained by filling in
the $c_{k}$-th box of the $1\times n$ diagram with $r_{k}$ for $1\leq k\leq
d $. Then a product of them will be denoted by a filling of a rectangular
diagram with multiple rows. We place the $i$-th factor in the $i$-th row
counting from the bottom row. For example, if $n=6$, then $%
[1,2,4,5|1,3,4,6]\times \lbrack 2,3,5,6|1,2,3,6]$ can be denoted by%
\begin{equation}  \label{tableau notation}
\tableau{ 2 &3 &5 & & &6\\1 & &2 & 4& &5}
\end{equation}

For $\mathbf{m}=(m_{1},\cdots ,m_{n-1})\in \mathbb{Z}_{>0}^{n-1}$, let $%
C^{(d)}=(1,2,\cdots ,d)$ and $R_{t,d}\subset \{1,\cdots ,n\}$ with $%
|R_{t,d}|=d$ for $1\leq d\leq n-1$. Then the product%
\begin{equation}
\prod_{1\leq t\leq m_{1}}[R_{t,1}:C^{(1)}]\times \cdots \times \prod_{1\leq
t\leq m_{n-1}}[R_{t,n-1}:C^{(n-1)}]  \label{standard monomials}
\end{equation}%
can be\ identified with a filling of $|\sum_{i}m_{i}|\times n$ rectangular
diagram.

With these fixed column indices $\left\{ C^{(d)}\right\} $, if the entries
in each row are strictly increasing from left to right and the entries in
each column is weakly increasing from top to bottom, then it can be
identified with a semistandard Young tableaux of shape $\lambda =(\lambda
_{1},\cdots ,\lambda _{n-1})$ with entries from $\{1,\cdots ,n\}$ in the
literature (e.g., \cite{Sta99}) where $\lambda _{i}=m_{i}+\cdots +m_{r}$ for 
$1\leq i\leq n-1$. Recall that they form a $\mathbb{C}$-basis of the section
space $H^{0}(G/B,\mathcal{O}_{\mathbf{m}})$, and called the \textit{standard
monomials} for $H^{0}(G/B,\mathcal{O}_{\mathbf{m}})$. See \cite{LG01}\cite%
{Se07} for further detail.

\medskip

\subsection{Schubert variety $X_{w}$\label{SchubertSMT}}

Moreover, standard monomial theory for $G/B$ descends to its Schubert
varieties $X_{w}$ in a way compatible with the embedding (\ref%
{flag-embedding}). Note that we can compare Pl\"{u}cker coordinates in the
Grassmannians in terms of the partial order given in (\ref{partial_order}).
Then for $X_{w}\subset G/B$, the kernel of the restriction map 
\begin{equation*}
H^{0}(G/B,\mathcal{O}_{\mathbf{m}})\rightarrow H^{0}(X_{w},\mathcal{O}_{%
\mathbf{m}})
\end{equation*}%
is spanned by $\prod_{d}\prod_{t}[R_{t,d}:C^{(d)}]$ such that $(w(1),\cdots
,w(d))\uparrow \nsucceq R_{t,d}$; and the following standard monomials in $%
H^{0}(G/B,\mathcal{O}_{\mathbf{m}})$ project to a $\mathbb{C}$-basis of $%
H^{0}(X_{w},\mathcal{O}_{\mathbf{m}})$:%
\begin{equation}
\left\{ \prod_{1\leq d\leq n-1}\prod_{1\leq t\leq
m_{d}}[R_{t,d}:C^{(d)}]:(w(1),\cdots ,w(d))\uparrow \succeq R_{t,d}\right\}
\label{Schubert condn}
\end{equation}%
See \cite{LG01} for further details.

\medskip

\subsection{Richardson variety $X_{w}^{v}$}

For $v\in \mathfrak{S}_{n}$, the opposite Schubert variety $X^{v}$ is the
Zariski closure of the $B^{-}$-orbit $B^{-}\cdot e_{v}$ where $B^{-}$ is the
opposite Borel subgroup. Then for $w,v\in \mathfrak{S}_{n}$, the Richardson
variety $X_{w}^{v}$ is defined as%
\begin{equation*}
X_{w}^{v}=X_{w}\cap X^{v}.
\end{equation*}%
In particular, this is non-empty if and only if $w\geq v$ with respect to
the Bruhat-Chevalley order.

By applying standard monomial theory to the opposite Schubert variety, we
can obtain standard monomials for $X_{w}^{v}$. This, via (\ref%
{flag-embedding}), follows directly from the case of the Richardson
varieties in the Grassmannian. For the restriction map%
\begin{equation*}
H^{0}(G/B,\mathcal{O}_{\mathbf{m}})\rightarrow H^{0}(X_{w}^{v},\mathcal{O}_{%
\mathbf{m}}),
\end{equation*}%
the following standard monomials in $H^{0}(G/B,\mathcal{O}_{\mathbf{m}})$
project to a $\mathbb{C}$-basis of $H^{0}(X_{w}^{v},\mathcal{O}_{\mathbf{m}%
}) $:%
\begin{equation}
\left\{ \prod_{1\leq d\leq n-1}\prod_{1\leq t\leq
m_{d}}[R_{t,d}:C^{(d)}]:(w(1),\cdots ,w(d))\succeq R_{t,d}\succeq
(v(1),\cdots ,v(d))\right\}  \label{Richardson condn}
\end{equation}%
See \cite{LG01} for further details.

\medskip

\section{Bott-Samelson varieties}

Let us review and generalize the results on the Bott-Samelson varieties
given in \cite{FK09}.

\subsection{Reduced word $\mathbf{i}$}

To obtain explicit descriptions, once and for all we fix the following
reduced decomposition of the longest element in $\mathfrak{S}_{n}$:%
\begin{equation*}
(s_{1})(s_{2}s_{1})\cdots (s_{n-1}s_{n-2}\cdots s_{1})
\end{equation*}%
where $s_{i}$ is the simple reflection $(i,i+1)$ for $1\leq i\leq n-1$, and
fix the following word%
\begin{equation}
\mathbf{i}=(i_{1},\cdots ,i_{\ell })  \label{reduced word}
\end{equation}%
associated with the above expression $s_{i_{1}}s_{i_{2}}\cdots s_{i_{\ell }}$
of the longest element. We write $w_{\mathbf{i}}$ for the longest element in 
$\mathfrak{S}_{n}$ and $\ell $ for the length of $w_{\mathbf{i}}$, which is $%
n(n-1)/2$.

\medskip

\subsection{Bott-Samelson varieties}

Let us consider a word $\mathbf{j}=(j_{1},\cdots ,j_{\ell ^{\prime }})$
whose corresponding expression of the element $w_{\mathbf{j}%
}=s_{j_{1}}s_{j_{2}}\cdots s_{j_{\ell ^{\prime }}}$ in $\mathfrak{S}_{n}$ is
reduced. The \textit{Bott-Samelson variety} is the quotient space%
\begin{equation*}
Z_{\mathbf{j}}=P_{j_{1}}\times P_{j_{2}}\times \cdots \times P_{j_{\ell
^{\prime }}}/B^{\ell ^{\prime }}
\end{equation*}%
where $P_{j_{i}}$ is the minimal parabolic subgroup of $G$ associated with
the simple reflection $s_{j_{i}}$, and $B^{\ell ^{\prime }}$ acts on the
product of $P_{j_{i}}$'s by 
\begin{equation*}
(p_{1},\cdots ,p_{\ell ^{\prime }}).(b_{1},\cdots ,b_{\ell ^{\prime
}})=(p_{1}b_{1},b_{1}^{-1}p_{2}b_{2},\cdots ,b_{\ell ^{\prime
}-1}^{-1}p_{\ell ^{\prime }}b_{\ell ^{\prime }}).
\end{equation*}

\medskip

We can also realize the Bott-Samelson variety $Z_{\mathbf{j}}$ as a
configuration variety \cite{Ma98}:%
\begin{equation}
Z_{\mathbf{j}}\subset Gr(\mathbf{j})=Gr(j_{1},n)\times \cdots \times
Gr(j_{\ell ^{\prime }},n)  \label{BSembedding}
\end{equation}%
More precisely, after realizing $Gr(j_{d},n)$ via $G/\widehat{P}_{j_{d}}$,
the Bott-Samelson variety $Z_{\mathbf{j}}$ is the closure of the $B$-orbit of%
\begin{equation*}
z_{\mathbf{j}}=(s_{j_{1}}\widehat{P}_{j_{1}},s_{j_{1}}s_{j_{2}}\widehat{P}%
_{j_{2}},...,s_{j_{1}}s_{j_{2}}\cdots s_{j_{\ell ^{\prime }}}\widehat{P}%
_{j_{\ell ^{\prime }}})
\end{equation*}%
Then for $\mathbf{m}=(m_{1},\cdots ,m_{\ell ^{\prime }})\in \mathbb{Z}%
_{>0}^{\ell ^{\prime }}$, we can consider a natural line bundle $L_{\mathbf{%
j,m}}$ induced from the Pl\"{u}cker bundles over $Gr(j_{i},n)$.

\medskip

\subsection{Homogeneous coordinate ring of $Z_{\mathbf{j}}$}

For any $\mathbf{m}=(m_{1},\cdots ,m_{\ell ^{\prime }})\in \mathbb{Z}%
_{>0}^{\ell ^{\prime }}$, with the realization of $Z_{\mathbf{j}}$ given in (%
\ref{BSembedding}), the section ring of $Z_{\mathbf{j}}$\ with respect to $%
L_{\mathbf{j,m}}$:%
\begin{equation*}
\mathcal{R}_{\mathbf{j,m}}=\bigoplus_{p\geq 0}H^{0}(Z_{\mathbf{j}},L_{%
\mathbf{j,m}}^{\otimes p})
\end{equation*}%
can be identified with a ring generated by products of determinants defined
by $\mathbf{j}$ and $\mathbf{m}$.

\begin{definition}
\label{ColumnSets}The column sets attached to $\mathbf{j}$ are%
\begin{equation*}
K_{\mathbf{j}}^{(r)}=s_{j_{1}}s_{j_{2}}\cdots s_{j_{r}}\{1,2,\cdots ,i_{r}\}
\end{equation*}%
for $1\leq r\leq \ell ^{\prime }$.
\end{definition}

With the notation we set in (\ref{RC_determinant}), let us consider a
multiset of determinants $[R_{t,r}:K_{\mathbf{j}}^{(r)}]$ whose column
indices are given by the column sets $K_{\mathbf{j}}^{(r)}$. Then, by
repeating $K_{\mathbf{j}}^{(r)}$'s $m_{r}$ times, the product $\mathcal{T}$
of determinants $[R_{t,r}:K_{\mathbf{j}}^{(r)}]$ for $1\leq t\leq m_{r}$ can
be encoded by a filling of a $|\mathbf{m|}\times n$\ rectangular diagram
having $[R_{t,r}:K_{\mathbf{j}}^{(r)}]$ as its $\left( m_{1}+\cdots
+m_{r-1}+t\right) $-th row counting from the bottom, where $|\mathbf{m|=}%
\sum_{i}m_{i}$.

\begin{definition}
A \textit{tableau }$\mathcal{T}$ of shape $(\mathbf{j,m)}$ is%
\begin{equation}
\mathcal{T}=\prod_{1\leq t\leq m_{1}}[R_{t,1}:K_{\mathbf{j}}^{(1)}]\cdot
\prod_{1\leq t\leq m_{2}}[R_{t,2}:K_{\mathbf{j}}^{(2)}]\cdot ...\cdot
\prod_{1\leq t\leq m_{\ell ^{\prime }}}[R_{t,\ell ^{\prime }}:K_{\mathbf{j}%
}^{(\ell ^{\prime })}]  \label{tableaux}
\end{equation}%
where for each $r$, all the row indexing sets satisfy $K_{\mathbf{j}%
}^{(r)}\succeq R_{t,r}$ for $1\leq t\leq m_{r}$. Let $\mathsf{M}(\mathbf{j,m}%
)$ be the space spanned by the tableaux of shape $(\mathbf{j,m)}$.
\end{definition}

As given in (\ref{tableau notation}), we will identify every tableau of
shape $(\mathbf{j,m)}$ with a filling of a rectangular diagram of size $|%
\mathbf{m|}\times n$. Hence the entry in the cell $(a,b)$ of a tableau $%
\mathcal{T}$ means the entry in the $a$-th row and $b$-th column in the
diagram realization of $\mathcal{T}$ counting from bottom to top and left to
right respectively.

\begin{proposition}[\S 3 \protect\cite{Ma98}]
\label{coordinate ring}For $\mathbf{m}=(m_{1},\cdots ,m_{\ell ^{\prime
}})\in \mathbb{Z}_{>0}^{\ell ^{\prime }}$, the section space $H^{0}(Z_{%
\mathbf{j}},L_{\mathbf{j,m}})$ of $Z_{\mathbf{j}}$ is isomorphic to the
space spanned by the tableaux of shape $(\mathbf{j,m)}$, i.e.,%
\begin{equation*}
H^{0}(Z_{\mathbf{j}},L_{\mathbf{j,m}})\cong \mathsf{M}(\mathbf{j,\mathbf{m}}%
).
\end{equation*}
\end{proposition}

\medskip

\subsection{Straight tableaux for $\mathcal{R}_{\mathbf{i,m}}$}

In the special case of $\mathbf{j=i}$, each column set $K_{\mathbf{i}}^{(r)}$
contains consecutive integers, and this fact lets us realize tableaux of
shape $(\mathbf{i,m)}$ in the context of tableaux of a \textit{row-convex
shape} which is defined in \cite{Ta01} as a generalized skew Young tableaux.
Using this observation, \cite{FK09} gives a presentation of the section ring%
\begin{equation*}
\mathcal{R}_{\mathbf{i,m}}=\bigoplus_{p\geq 0}\mathsf{M}(\mathbf{i,}p\mathbf{%
m})
\end{equation*}%
in terms of tableaux, and then identified an explicit basis. Note that up to
sign, we can always assume that the entries in each row of $\mathcal{T}$ are
increasing from left to right. If such is the case, then $\mathcal{T}$ is
called a \textit{row-standard tableau}.

\begin{definition}
A row-standard \textit{tableau} $\mathcal{T}$ of shape $(\mathbf{i,m)}$ is
called a straight tableau, if $\mathcal{T}$ as a $|\mathbf{m}|\times n$
tableau (\ref{tableaux}) satisfies the following condition: for two cells $%
(i,k)$ and $(j,k)$ with $i<j$ in the same column, the entry in the cell $%
(i,k)$ may be strictly larger than the entry in $(j,k)$ only if the cell $%
(i,k-1)$ exists and contains an entry weakly larger than the one in the cell 
$(j,k)$.
\end{definition}

\begin{theorem}[\protect\cite{FK09}]
\label{STM_BSi}Straight tableaux form a $\mathbb{C}$-basis for the $\mathbb{Z%
}$-graded algebra $\mathcal{R}_{\mathbf{i,m}}$. In particular, straight
tableaux of shape $(\mathbf{i,}p\mathbf{\mathbf{m})}$ form a $\mathbb{C}$%
-basis for the section space $\mathsf{M}(\mathbf{i,}p\mathbf{\mathbf{m}})$.
\end{theorem}

Our next task is to extend the above result to $\mathcal{R}_{\mathbf{j,m}}$
for a subword $\mathbf{j}$ of $\mathbf{i}$.

\medskip

\subsection{Relative description of $\mathcal{R}_{\mathbf{j,m}}$ to $%
\mathcal{R}_{\mathbf{i,m}}$}

For a subword $\mathbf{j}$ of $\mathbf{i}$, we can study a relative
description of $\mathcal{R}_{\mathbf{j,m}}$ to $\mathcal{R}_{\mathbf{i,m}}$
by using the canonical embedding $Z_{\mathbf{j}}\subset Z_{\mathbf{i}}:$%
\begin{diagram}
Gr(\mathbf{j})  &      \rEmbed   &  Gr(\mathbf{i}) \\
\uEmbed          &               &  \uEmbed       \\
Z_{\mathbf{j}}  &      \rEmbed   &  Z_{\mathbf{i}}
\end{diagram}%
See (\ref{BSembedding}) for notation.

For this purpose, in what follows, we write a subword $\mathbf{j}$ of $%
\mathbf{i}$ as $(j_{1},\cdots ,j_{\ell })$ by adopting the convention of
using $0$ for the omitted letters. Its associated element $w_{\mathbf{j}}$
in $\mathfrak{S}_{n}$ is $s_{j_{1}}s_{j_{2}}\cdots s_{j_{\ell }}$ with $s_{0}
$ being the identity in $\mathfrak{S}_{n}$. We further assume that $\mathbf{j%
}$ is reduced. For the rest of our discussion, if $j_{r}=0$, then we assume
the corresponding object is considered to be omitted or a trivial one. For
example, in a product of the Grassmannians $Gr(j_{1},n)\times
Gr(j_{2},n)\times \cdots \times Gr(j_{\ell },n)$, if $j_{r}=0$ then we omit $%
Gr(j_{r},n)$ and therefore the column set $K_{\mathbf{j}}^{(r)}=\emptyset $.
Also, in $\mathbf{m}=(m_{1},\cdots ,m_{\ell })$ attached to $\mathbf{j}$, we
assume $m_{r}=0$ if $j_{r}=0$.

\medskip

Let us consider the restriction map $H^{0}(Z_{\mathbf{i}},L_{\mathbf{i,m}%
})\rightarrow H^{0}(Z_{\mathbf{j}},L_{\mathbf{j,m}})$, or more explicitly
the following map $\phi :\mathsf{M}(\mathbf{i,m})\rightarrow \mathsf{M}(%
\mathbf{j,m})$%
\begin{eqnarray}
&&\phi \left( \prod_{1\leq t\leq m_{1}}[R_{t,1}:K_{\mathbf{i}}^{(1)}]\cdot
...\cdot \prod_{1\leq t\leq m_{\ell }}[R_{t,\ell }:K_{\mathbf{i}}^{(\ell
)}]\right)  \label{restriction map} \\
&=&\prod_{1\leq t\leq m_{1}}\phi _{1}([R_{t,1}:K_{\mathbf{i}}^{(1)}])\cdot
...\cdot \prod_{1\leq t\leq m_{\ell }}\phi _{\ell }([R_{t,\ell }:K_{\mathbf{i%
}}^{(\ell )}])  \notag
\end{eqnarray}%
where $\phi _{r}([R_{t,r}:K_{\mathbf{i}}^{(r)}])=[R_{t,r}:K_{\mathbf{j}%
}^{(r)}]$ for $r$ such that $j_{r}=i_{r}$; and $\phi _{r}([R_{t,r}:K_{%
\mathbf{i}}^{(r)}])=1$ for $r$ such that $j_{r}=0$. Then $\phi $ is
surjective and its kernel is the union of the kernels of $\phi _{r}$'s.

\begin{proposition}
\label{SMT_Zj}i) The kernel of the restriction map $\phi :H^{0}(Z_{\mathbf{i}%
},L_{\mathbf{i,m}})\rightarrow H^{0}(Z_{\mathbf{j}},L_{\mathbf{j,m}})$ is
spanned by%
\begin{equation*}
\left\{ \prod_{1\leq t\leq m_{1}}[R_{t,1}:K_{\mathbf{i}}^{(1)}]\cdot
...\cdot \prod_{1\leq t\leq m_{\ell }}[R_{t,\ell }:K_{\mathbf{i}}^{(\ell
)}]\in \mathsf{M}(\mathbf{i,m}):R_{t,r}\npreceq K_{\mathbf{j}}^{(r)}\text{
for }r\text{ such that }j_{r}=i_{r}\right\}
\end{equation*}%
ii) The straight tableaux of shape $(\mathbf{i,m})$ such that $%
R_{t,r}\preceq K_{\mathbf{j}}^{(r)}$ for $r$ such that $j_{r}=i_{r}$ project
to a $\mathbb{C}$-basis of the space $H^{0}(Z_{\mathbf{j}},L_{\mathbf{j,m}})$%
.
\end{proposition}

\begin{proof}
From the embedding $Z_{\mathbf{j}}\subset Gr(\mathbf{j})$ in (\ref%
{BSembedding}), it is enough to check the statements for the individual
factors of $Gr(\mathbf{j})$. To each factor $G/\widehat{P}_{j_{d}}$ of $Gr(%
\mathbf{j})$, the Bott-Samelson variety $Z_{\mathbf{j}}$ projects to a
Schubert variety, because it is the closure of the $B$ orbit of the $T$%
-invariant element $s_{j_{1}}s_{j_{2}}\cdots s_{j_{r}}\widehat{P}_{j_{r}}$
in $G/\widehat{P}_{j_{r}}$. Then the statements follow from the standard
monomial description for the section space of Schubert varieties given in
\S \ref{SchubertSMT} and Theorem \ref{STM_BSi}.
\end{proof}

The map $\phi $ naturally extends to $\Phi $ from $\mathcal{R}_{\mathbf{i,m}%
}=\bigoplus_{p\geq 0}\mathsf{M}(\mathbf{i,}p\mathbf{m})$ to $\mathcal{R}_{%
\mathbf{j,m}}$. Therefore, the above description gives the relative
description of $\mathcal{R}_{\mathbf{j,m}}$ to $\mathcal{R}_{\mathbf{i,m}}$.
From now on, we identify the space $\mathsf{M}(\mathbf{j,m})$ with the
quotient space 
\begin{equation*}
\mathsf{M}(\mathbf{j,m})=\mathsf{M}(\mathbf{i,m})/\ker \phi .
\end{equation*}

\begin{example}
\label{row convex shape}For $n=4,$ $\mathbf{i}=(1,2,1,3,2,1)$ and the column
sets $K_{\mathbf{i}}^{(r)}$ are 
\begin{eqnarray*}
&&K_{\mathbf{i}}^{(1)}=\{2\},K_{\mathbf{i}}^{(2)}=\{2,3\},K_{\mathbf{i}%
}^{(3)}=\{3\}, \\
&&K_{\mathbf{i}}^{(4)}=\{2,3,4\},K_{\mathbf{i}}^{(5)}=\{3,4\},K_{\mathbf{i}%
}^{(6)}=\{4\}.
\end{eqnarray*}%
For the subword $\mathbf{j}=(1,0,0,0,2,1)$, its column sets $K_{\mathbf{j}%
}^{(r)}$ are%
\begin{eqnarray*}
&&K_{\mathbf{j}}^{(1)}=\{2\},K_{\mathbf{j}}^{(2)}=\emptyset ,K_{\mathbf{j}%
}^{(3)}=\emptyset , \\
&&K_{\mathbf{j}}^{(4)}=\emptyset ,K_{\mathbf{j}}^{(5)}=\{2,3\},K_{\mathbf{j}%
}^{(6)}=\{3\}.
\end{eqnarray*}%
With $\mathbf{m}=(1,1,1,1,1,1)$, the surjection $\Phi $ sends tableaux shape
of $(\mathbf{i,m})$\ to tableaux of shape $(\mathbf{j,m})$ as follows:%
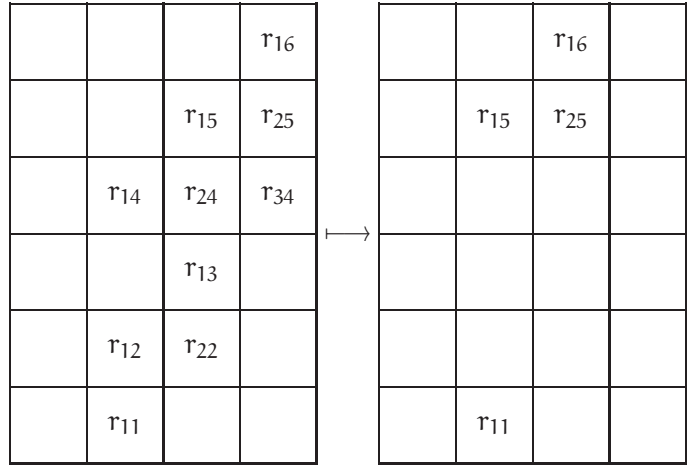
\begin{equation*}
\tableau{ & & & r_{16} \\ & & r_{15} & r_{25} \\ & r_{14} & r_{24} & r_{34}
\\ & & r_{13} & \\ & r_{12} & r_{22} & \\ & r_{11} & & }\longmapsto %
\tableau{ & & r_{16} & \\ & r_{15} & r_{25}& \\ & & & \\ & & & \\ & & & \\ &
r_{11} & & }
\end{equation*}%
Then the kernel is spanned by the tableaux with $(r_{15},r_{25})\npreceq
(2,3)$ and $(r_{11})\npreceq (2)$. The second condition in this case is void.
\end{example}

\medskip

\section{RR Varieties}

\subsection{Involution $w_{\mathbf{i}}$}

Let us state parallel results for opposite Schubert varieties and the
corresponding Bott-Samelson varieties. Fix a subword $\mathbf{j}%
=(j_{1},\cdots ,j_{\ell })$ of $\mathbf{i}$ such that the corresponding
expression of the element $w_{\mathbf{j}}=s_{j_{1}}s_{j_{2}}\cdots
s_{j_{\ell }}$ in $\mathfrak{S}_{n}$ is reduced. Then we define the
corresponding opposite Schubert variety as%
\begin{equation*}
X^{\mathbf{j}}=\overline{B^{-}\cdot e_{v}}
\end{equation*}%
where $v=w_{\mathbf{i}}w_{\mathbf{j}}$, and its corresponding Bott-Samelson
variety $Z^{\mathbf{j}}$ as the closure of $B^{-}$ orbit in the product of
the Grassmannians $Gr(j_{r},n)$ as in (\ref{BSembedding}).

The section space of the line bundle $\widetilde{L}_{\mathbf{j,m}}$ over $Z^{%
\mathbf{j}}$ can be obtained by applying the involution $w_{\mathbf{i}}$ to
the line bundle $L_{\mathbf{j,m}}$ over $Z_{\mathbf{j}}.$ A row-standard
tableau $\mathcal{T}$ of the form (\ref{tableaux}) is sent to%
\begin{equation*}
w_{\mathbf{i}}(\mathcal{T)}=\left( \prod_{m_{\ell }\geq t\geq 1}[\widetilde{R%
}_{t,\ell }:\widetilde{K}_{\mathbf{j}}^{(\ell )}]\right) \cdot \left(
\prod_{m_{\ell -1}\geq t\geq 1}[\widetilde{R}_{t,\ell -1}:\widetilde{K}_{%
\mathbf{j}}^{(\ell -1)}]\right) \cdot ...\cdot \left( \prod_{m_{1}\geq t\geq
1}[\widetilde{R}_{t,1}:\widetilde{K}_{\mathbf{j}}^{(1)}]\right) 
\end{equation*}%
where for a subset $X=\{x_{1}<\cdots <x_{r}\}$ of $\{1,\cdots ,n\}$, 
\begin{equation*}
\widetilde{X}=\{n+1-x_{r},\cdots ,n+1-x_{1}\}.
\end{equation*}%
\ Note that we reverse the order of multiplication. In terms of our diagram
notation (\ref{tableau notation}), $w_{\mathbf{i}}(\mathcal{T)}$ is obtained
by rotating $\mathcal{T}$ by $180%
{{}^\circ}%
$ and then replacing row indexing entries $r_{ij}$ by $n+1-r_{ij}$.

\begin{example}
The surjection in Example \ref{row convex shape}, via the involution $w_{%
\mathbf{i}}$, corresponds the restriction map from $H^{0}(Z^{\mathbf{i}},%
\widetilde{L}_{i\mathbf{,m}})$ to $H^{0}(Z^{\mathbf{j}},\widetilde{L}_{%
\mathbf{j,m}})$:%
\begin{equation*}
\tableau{ & & r'_{11} & \\ & r'_{22} & r'_{12} & \\ & r'_{13} & & \\ r'_{34}
& r'_{24} & r'_{14} & \\ r'_{25} & r'_{15} & & \\ r'_{16} & & & }\longmapsto %
\tableau{ & & r'_{11} & \\ & & & \\ & & & \\ & & & \\ r'_{25} & r'_{15} & &
\\ r'_{16} & & & }
\end{equation*}%
where $r_{ij}^{\prime }=n+1-r_{ij}$.
\end{example}

Then it follows from Proposition \ref{SMT_Zj} that

\begin{corollary}
For $\mathbf{m}=(m_{1},\cdots ,m_{\ell })\in \mathbb{Z}_{>0}^{\ell }$
attached to $\mathbf{j}$, we have%
\begin{equation*}
H^{0}(Z^{\mathbf{j}},\widetilde{L}_{\mathbf{j,m}})\cong w_{\mathbf{i}}(%
\mathsf{M}(\mathbf{j,m})),
\end{equation*}%
and for the straight tableaux $\mathcal{T}$ of shape $(\mathbf{j,m})$, $w_{%
\mathbf{i}}(\mathcal{T)}$\ project to a $\mathbb{C}$-basis of $H^{0}(Z^{%
\mathbf{j}},\widetilde{L}_{\mathbf{j,m}})$.
\end{corollary}

\medskip 

\subsection{Fiber Product}

For $\mathbf{m}=(m_{1},\cdots ,m_{\ell })\in \mathbb{Z}_{>0}^{\ell }$, let
us consider the homogeneous coordinate rings $\mathcal{R}_{\mathbf{j,m}}$
and $\mathcal{R}^{\mathbf{k,m}}$ of $Z_{\mathbf{j}}$ and $Z^{\mathbf{k}}$
respectively%
\begin{eqnarray*}
\mathcal{R}_{\mathbf{j,m}} &=&\bigoplus_{p\geq 0}\mathsf{M}(\mathbf{j,}p%
\mathbf{m}) \\
\mathcal{R}^{\mathbf{k,m}} &=&\bigoplus_{p\geq 0}w_{\mathbf{i}}(\mathsf{M}(%
\mathbf{k,}p\mathbf{m}))
\end{eqnarray*}

Then by taking the last $(n-1)$ entries of $\mathbf{m}$, we set $\mathbf{m}%
_{0}=(m_{\ell -n+2},\cdots ,m_{\ell })$. Then let $\mathcal{A}=\mathcal{A}%
_{w}^{v}=\bigoplus_{p\geq 0}H^{0}(X_{w}^{v},\mathcal{O}_{\mathbf{m}%
_{0}}^{\otimes p})$ be the homogeneous coordinate ring, given in (\ref%
{Richardson condn}), of the Richardson variety $X_{w}^{v}=X_{w}\cap X^{v}$
where $w=w_{\mathbf{j}}$ and $v=w_{\mathbf{i}}w_{\mathbf{k}}$.

We define the coproduct of $\mathcal{R}_{\mathbf{j,m}}$ and $\mathcal{R}^{%
\mathbf{k,m}}$ over $\mathcal{A}$%
\begin{equation*}
\mathcal{R}_{\mathbf{j}}^{\mathbf{k}}=\mathcal{R}_{\mathbf{j,m}}\otimes _{%
\mathcal{A}}\mathcal{R}^{\mathbf{k,m}}
\end{equation*}%
with respect to the following injective maps: $\varphi _{\mathbf{j}%
}:H^{0}(X_{w}^{v},\mathcal{O}_{\mathbf{m}_{0}})\rightarrow \mathsf{M}(%
\mathbf{j,m})$ sending $\mathcal{T}$ to

\begin{equation*}
\varphi _{\mathbf{j}}\left( \mathcal{T}\right) =\left( \prod_{1\leq t\leq
m_{1}}[R_{t,1}^{0}:K_{\mathbf{j}}^{(1)}]\cdot ...\cdot \prod_{1\leq t\leq
m_{\ell +1-n}}[R_{t,\ell +1-n}^{0}:K_{\mathbf{j}}^{(\ell +1-n)}]\right)
\cdot \mathcal{T}
\end{equation*}%
in the quotient $\mathsf{M}(\mathbf{j,m})=\mathsf{M}(\mathbf{i,m})/\ker \phi 
$ where $R_{t,r}^{0}=\{1,2,\cdots ,|K_{\mathbf{j}}^{(r)}|\}$ for $1\leq
r\leq \ell +1-n$, and 
\begin{equation*}
\varphi ^{\mathbf{k}}:H^{0}(X_{w}^{v},\mathcal{O}_{\mathbf{m}%
_{0}})\rightarrow w_{\mathbf{i}}(\mathsf{M}(\mathbf{k,m}))
\end{equation*}%
defined by $\varphi ^{\mathbf{k}}(\mathcal{T})=w_{\mathbf{i}}(\varphi _{%
\mathbf{j}}(\mathcal{T}))$. Note that these maps correspond to the
projections from the Bott-Samelson varieties to the flag varieties (cf. \cite%
[\S 4.3]{FK09}).

With $(\mathcal{R}_{\mathbf{j,m}},\mathcal{R}^{\mathbf{k,m}},\mathcal{A}%
_{w}^{v},\varphi _{\mathbf{j}},\varphi ^{\mathbf{k}})$, we define the 
\textit{RR variety} as the fiber product of $Z_{\mathbf{j}}$ and $Z^{\mathbf{%
k}}$ over $X_{w}^{v}$.

\medskip 

\subsection{Toric Degenerations}

A monomial order on the polynomial ring $\mathbb{C}[M_{n}]$ is called a 
\textit{diagonal term order} if the leading monomial of a determinant of any
minor over $M_{n}$ is equal to the product of the diagonal elements. For a
subring $\mathcal{R}$ of the polynomial ring we let $in(\mathcal{R})$ denote
the algebra generated by the leading monomials $in(f)$ of all $f\in \mathcal{%
R}$ with respect to a given monomial order. The leading monomials for all $%
f\in \mathcal{R}$ form an affine semigroup, therefore $in(\mathcal{R})$ is a
semigroup algebra representing a toric variety in the sense of \cite{St95}.

In \cite{Ta01}, it is shown that for a row-convex shape $(\mathbf{h,m)}$,
straight tableaux of shape $(\mathbf{h,m)}$ form a SAGBI basis of the graded
algebra $\mathcal{R}\subset \mathbb{C}[M_{n}]$ generated by tableaux of
shape $(\mathbf{h,m)}$ with respect to any diagonal term order. A finite
SAGBI basis for $\mathcal{R}$ provides a toric degeneration of $Spec(%
\mathcal{R})$ to $Spec(in(\mathcal{R}))$ or ${Proj}(\mathcal{R})$ to 
${Proj}(in(\mathcal{R}))$ if there is a $\mathbb{Z}$-grading. This is called a 
\textit{SAGBI-degeneration} (cf. \cite[p.281]{MS05}\cite[Theorem 1.2]{CHV96}%
). Using this method, \cite{FK09} shows $Z_{\mathbf{i}}$ can be flatly
deformed into a toric variety. For $Z_{\mathbf{j}}$, we will apply an
analogous method to the quotient algebra $\mathcal{R}_{\mathbf{j,m}}=%
\mathcal{R}_{\mathbf{i,m}}/\ker \Phi $.

\begin{theorem}
\label{deformation}The Bott-Samelson variety $Z_{\mathbf{j}}$ can be flatly
deformed into a toric variety.
\end{theorem}

\begin{proof}
We show that there is a flat $\mathbb{C}[t]$ module $\mathcal{R}_{\mathbf{j,m%
}}^{t}$ whose general fiber is isomorphic to $\mathcal{R}_{\mathbf{j,m}}$
and special fiber is isomorphic to a semigroup ring corresponding to an
initial object of $\mathcal{R}_{\mathbf{j,m}}$. To specify this initial
object of the quotient algebra, we will use the fact that every element of $%
\mathcal{R}_{\mathbf{j,m}}$ has a cannonical representative. From
Proposition \ref{SMT_Zj}, every homogeneous element $H$ of the quotient $%
\mathcal{R}_{\mathbf{i,m}}/\ker \Phi $ can be expressed as a linear
combination of straight tableaux of the same shape%
\begin{equation*}
H=\sum c_{i}\mathcal{T}_{i}
\end{equation*}%
From the fact that leading monomials of straight tableaux of a fixed shape
are distinct (\cite{Ta01}), the leading monomial of $f$ should be equal to
the leading monomial of $\mathcal{T}_{i}$ for some $i$. Therefore, we have a
well defined notion of the leading monomials $in(H)$ for $H$ in $\mathcal{R}%
_{\mathbf{i,m}}/\ker \Phi $ and it is equal to $in(\mathcal{T})$ for a
straight tableau $\mathcal{T}$. In this sense, straight tableaux form a
SAGBI basis for $\mathcal{R}_{\mathbf{i,m}}/\ker \Phi $, i.e., the semigroup
ring $in(\mathcal{R}_{\mathbf{i,m}}/\ker \Phi )$ is generated by the leading
monomials of straight tableaux. Since for any shape $(\mathbf{j,m)}$ we only
have a finite number of straight tableaux, straight tableaux form a finite
SAGBI\ basis for $\mathcal{R}_{\mathbf{i,m}}/\ker \Phi $. Then, by the same
argument given in \cite[Proposition 1]{ST99}, we have a $\mathbb{Z}_{\geq 0}$%
-filtration $\{F_{\alpha }\}$ on $\mathcal{R}_{\mathbf{j,m}}\cong \mathcal{R}%
_{\mathbf{i,m}}/\ker \Phi $ such that the Rees algebra $\mathcal{R}_{\mathbf{%
j,m}}^{t}$ of $\mathcal{R}_{\mathbf{j,m}}$ with respect to $\{F_{\alpha }\}$:%
\begin{equation*}
\mathcal{R}_{\mathbf{j,m}}^{t}=\bigoplus_{\alpha \geq 0}F_{\alpha }(\mathcal{%
R}_{\mathbf{j,m}})t^{\alpha }
\end{equation*}%
which is flat over $\mathbb{C}[t]$, has a general fiber isomorphic to $%
\mathcal{R}_{\mathbf{j,m}}$ and the special fiber isomorphic to the
semigroup ring $in(\mathcal{R}_{\mathbf{i,m}}/\ker \Phi )$.
\end{proof}

\begin{corollary}
The RR variety can be flatly deformed into a fiber product of toric
varieties.
\end{corollary}

\end{document}